\documentclass[a4paper,12pt]{article}
  
\usepackage{amsfonts}
\usepackage{amssymb,amsmath} 
\numberwithin{equation}{section}
\newtheorem{theorem}{Theorem}[section]

\newtheorem{proposition}{Proposition}[section]

\usepackage{cite}
\usepackage{hyperref}
\hypersetup{urlcolor=blue, colorlinks=true, citecolor=blue, linkcolor=blue}

\begin{document}

\begin{center}
\textbf{\Large{Constant Curvature Conditions For Generalized Kropina Spaces}}\\
\end{center}

\begin{flushleft}
\textbf{Gauree Shanker$^1$, Ruchi Kaushik Sharma$^2$}\\

\scriptsize{$^1$ Department of Mathematics and Statistics, Central University of Punjab, Bhatinda-151001, India; Email:grshnkr2007@gmail.com}\\
\scriptsize{$^2$ Department of Mathematics \& Statistics, Banasthali University, Banasthali 304 022, Rajasthan, India}; Email:ruchikaushik07@gmail.com.\\
[0.2cm]
\end{flushleft}


\textbf{Abstract -}The classification of Finsler spaces of constant curvature is an interesting and important topic of research in differential geometry. In this paper we obtain necessary and sufficient conditions for generalized Kropina space to be of constant flag curvature.  \\
\\
\textbf{M. S. C. 2010}: 58B20, 53B21, 53B40, 53C60.\\
{\bf Keywords -} Riemannian spaces; Killing vector fields; Finsler metrics; Kropina metrics; Generalized Kropina metrics.

\section{Introduction}
Finsler metrics are generalization of Riemannian metrics in the sense that they depend both on the position and direction while its counterpart depend only on position. Generalized Kropina metric belongs to the large class of ($\alpha, \beta$)-metric, first introduced by Matsumoto \cite{Mats92}. ($\alpha, \beta$)-metrics are constructed by Riemannian metric $\alpha = \sqrt{a_{ij}(x)y^iy^j}$ and the differential 1-form $\beta = b_i(x)y^i$. Some notable ($\alpha, \beta$)-metrices are: Randers metric: F = $\alpha$ + $\beta$; Kropina metric: F = $\frac{\alpha^2}{\beta}$; generalized Kropina metric: F = $\frac{\alpha^{m+1}}{\beta^m}$ (m $\neq$ 0. -1); Matsumoto metric : F = $\frac{\alpha^2}{\alpha - \beta}$ and square metric: F = $\frac{(\alpha + \beta)^2}{\alpha}$. Contrary to other ($\alpha$, $\beta$)-metrics, Kropina metric and generelised Kropina metric are not regular but they have wide applications in other branches of science. \\
Classification of ($\alpha, \beta$)-metrics is very important problem. Several geometers have worked on this topic with different perspective.  

\section{The description of generalized Kropina metric}
Let (M, $\alpha$) be an n($\geq$ 2)-dimensional smooth manifold endowed with Riemannian metric $\alpha$. A generalized Kropina space $\Big(M, \dfrac{\alpha^{m+1}}{\beta^m}\Big)$ is a Finsler space whose fundamental function is given by F = $\frac{\alpha^{m+1}}{\beta^m}$ ( m $\neq$ 0, -1), where $\alpha = \sqrt{a_{ij}(x)y^iy^j}$ is a Riemannian metric and $\beta = b_i(x)y^i$ is a differential 1-form. For our purpose we assume that the matrix $(a_{ij})$ is positive definite.
\\
It is to be noted that Randers spaces (M, F = $\alpha + \beta$) on TM are Lagragian duals of Kropina spaces ($M = \bar{F} = \frac{\bar{\alpha}^2}{\beta}$) on $T^{*}M$ and vice-versa in the case $b^2$ = 1, where b is the Riemannian length of $\beta$. \par
L-dual for generalized Kropina spaces and some other Finsler spaces have been obtained in (\cite{Shankr11}, \cite{Shanker15}, \cite{ShanRavi11}). Furthermore, for regular Lagrangians, the necessary and sufficient  condition for a Finsler space to be of constant flag curvature K is that its dual space is also of constant flag curvature $\bar{K}$ (\cite{Shimada96}, \cite{Shimada97}).
\par
Define a new Reimannian metric h = $\sqrt{h_{ij}(x)y^iy^j}$ and a vector field W = $W^i(\frac{\partial}{\partial x^i})$ on M by\\

$h_{ij} = e^{k(x)}a_{ij}$ and $W_i = e^{k(x)}b_i$,\\

where
 $W_i = h_{ij}W^j$.
\\
Then the generalized Kropina metric can be written as
\begin{align}
F = \dfrac{\alpha^{m+1}}{\beta^m} = \pi\frac{h_{00}^{\frac{m+1}{2}}}{W_{0}^m},
\end{align}
where $\pi = \frac{e^{\frac{(m-1)k(x)}{2}}}{2^m}$.
\section{The coefficients of the geodesic spray}
Let us recall \cite{yoshi12}, the following theorem for later use: 
\begin{theorem}
	Let (M, g) and (M, $g^* = e^{\rho}$g), where g = $\sqrt{g_{ij}(x)y^iy^j}$ and $g^* = \sqrt{g_{ij}^*(x)y^iy^j}$ respectively, be two n-dimensional Riemannian spaces which are conformal to each other. Furthermore, let ${{\gamma_{j}}^{ i}}_k$ and ${{\gamma_{j}}^{*i}}_k$ be the coefficients of Levi-Civita connection of (M, g) and (M, $g^*$), respectively. Then, we have 
	\begin{align}
	{g^*}_{ij} = e^{2\rho}g_{ij}, g^{*ij} = e^{-2\rho}g^{ij}
	\end{align}
	and
	\begin{align} 
	{{\gamma_{j}}^{*i}}_k = {{\gamma_{j}}^{i}}_k + \rho_j\delta^i_k + \rho_k\delta^i_j - \rho^ig_{jk},
	\end{align}
	where $\rho_i = \frac{\partial \rho}{\partial x^i}$
	$\rho^i = g^{ij}\rho_j$.
\end{theorem}
From (2.1), we have $h_{ij} = e^ka_{ij}$. Applying Theorem 3.1, we get
\begin{align}
^{h}{{\gamma_{j}}^{i}}_k =  ^{\alpha}{{\gamma_{j}}^{i}}_k + \frac{1}{2}k_j\delta^i_k + \frac{1}{2}k_k\delta^i_j - \frac{1}{2}k^ia_{jk},
\end{align}
where $^{h}{{\gamma_{j}}^{i}}_k \quad \text{and} \quad ^{\alpha}{{\gamma_{j}}^{i}}_k$ are the coefficients of levi-Civita connection of (M, h) and (M, $\alpha$) respectively,  $k_i = \partial k/\partial x^i$  and $k^i = a^{ij}k_j$. Transvecting (3.1) by $y^jy^k$, we get
\begin{align}
^{h}{\gamma_{0}^{i}}_0 =  ^{\alpha}{\gamma_{0}^{i}}_0 + k_0y^i - \frac{1}{2}h_{00}\bar{k}^i,
\end{align}
where $\bar{k}^i = h^{ij}k_j$ and the index 0 means the transvection by $y^i$. 
\par
we denote the covariant derivative in the Riemannian space (M, $\alpha$) by (;i) and introduce the following notations: $s_{ij} := \frac{b_{i;j} - b_{j;i}}{2}$, $r_{ij} = \frac{b_{i;j} + b_{j;i}}{2}$, $s_j := b^i s_{ij} $.
\par
In \cite{Basco07}, the authors have shown that the coefficients $G^i$ of geodesic spray in a Finsler space (M, F = $\alpha \phi(s)$), where s = $\beta/\alpha$ and $\phi$ is a differential function of s alone, are given by
\begin{align}
2G^i = ^{\alpha}{\gamma_{0}^{i}}_0 + 2\omega \alpha s_0^i + 2\Theta(r_{00} - 2\alpha \omega s_0) \Bigl(\frac{y^i}{\alpha} + \frac{\omega'}{\omega - s\omega'}b^i\Bigr),
\end{align}

where $\omega := \frac{\phi'}{\phi - s\phi'}$,
and $\Theta := \frac{\omega - s\omega'}{2{1 + s\omega + (b^2 - s^2)\omega'}}$.
\begin{align}
2G^i = ^{h}{\gamma_{0}^{i}}_0 - k_0y^i + \frac{1}{2}h_{00}\bar{k}^i - F s_0^i - \frac{1}{b^2}(r_{00} + Fs_0)\Bigl(\frac{2}{F}y^i - b^i\Bigr).
\end{align}
For a Generalized Kropina space $\bigl(M, \frac{\alpha^{m+1}}{\beta^m}\bigr)$ a new Riemannian metric  h = $\sqrt{h_{ij}y^iy^j}$ and a vector field W = $W^i(\partial / \partial x^i)$ are defined by (2.1).
\\
So, the vector field W satisfies the condition $\|W\|$ = 1 and we have F = $\pi\frac{h_{00}^{\frac{m+1}{2}}}{W_0^m}$.
\par
Therefore, we get

\begin{align}
2G^i = ^{h}{\gamma_{0}^{i}}_0 + 2\Phi^i,
\end{align}
where
\begin{align}
2\Phi^i : &= -k_0 y^i + \frac{1}{2}h_{00}\bar{k}^i - \dfrac{2m\pi}{m+1}\frac{h_{00}^{\frac{m+1}{2}}}{W_0^m} s_0^i 
\nonumber \\
&
- \frac{2ms}{s^2 - ms^2 + mb^2} \Bigl(r_{00} + \dfrac{2m \pi}{m+1}\frac{h_{00}^{\frac{m+1}{2}}}{W_0^m} s_0\Bigr),
\end{align}

\begin{align*}
\phi(s) =& \frac{1}{s^m},\\
\omega  =& -\frac{m}{1+m}\frac{1}{s},\\
\omega' =& \frac{m}{1+m}\frac{1}{s^2},\\
\Theta  =& -\frac{ms}{s^2 - ms^2 +mb^2}.
\end{align*} 
\textbf{Remark.}
	We can introduce a Finsler connection $\Gamma^* = (^{h}{\gamma_{j}^{i}}_k(x), N^{i}_{j} := {^{h}\gamma_{j}^{i}}_k(x)y^k, {C_{j}^ { i}}_ {k}$) associated with the linear connection ${^{h}\gamma_{j}^{i}}_k(x)$ of the Riemannian space (M, h).\\
	The h-covariant derivative are defined as follows \cite{Mats86}:
	\par
	For a vector field $W^i(x)$ on M, \\
	(1) $W^i(x)_{\rVert j} = \frac{\partial W^i}{\partial x^j} - \frac{\partial W^i}{\partial y^s} N^s_j + {^{h}\gamma_{j}^{i}}_s W^s = \frac{\partial W^i}{\partial x^j} + {^{h}\gamma_{j}^{i}}_s W^s.$\\
	For a reference vector $y^i$,\\
	(2) $y^i_{\rVert j} = \frac{\partial y^i}{\partial x^j} - \frac{\partial y^i}{\partial y^s}N_j^s + {^{h}\gamma_{j}^{i}}_s y^s = -N_i^j + N_i^j = 0.$\\
	
	We put 
	\begin{align*}
	R_{ij} := \frac{W_{i\rVert j} + W_{j \rVert i}}{2}, \quad S_{ij} := \frac{W_{i\rVert j} - W_{j \rVert i}}{2}, \quad R^i_j := h^{ir}R_{rj}, \quad S_i^j =h^{ir}S_{rj},\\
	R_i := W^r R_{ri}, \quad S_i := W^rS_{ri}, \quad R^i := h^{ir}R_r, \quad S^i := h^{ir}S_r.\\
	\text{It follows} \quad r_{ij} = 2e^{-k}\bigl(R_{ij} - \frac{1}{2}W_r\bar{k}^rh_{ij}\bigr), \quad s_{ij} = 2e^{-k}\bigl(S_{ij} + \frac{ k_iW_j - k_jW_i}{2}\bigr).\\
	\end{align*}
	Furthermore, we get
	\begin{align*}
	& s_j^i = 2S_i^j + \bar{k}^i W_j - k_jW^i, \nonumber \\
	 & s_0^i = 2S_0^i + W_0\bar{k}^i - k_W^i,
	\nonumber
	\\
	& s^0_i = 2e^{-k}(2S_i + W_r\bar{k}^rW_i - k_i), \nonumber \\ 
	&s_0 = 2e^{-k}(2S_0 + W_r\bar{k}^rW_0 - k_0), \nonumber \\
	& r_{00} = 2e^{-k}\bigl(R_{00} - \frac{1}{2}W_r\bar{k}^rh_{00}\bigr), \nonumber \\
	 &\quad b^i = a^{ir}b_r = e^kh^{ir}\frac{2W_r}{e^k} = 2W^i.
	\end{align*} 
	Substituting above equalities in (3.8), we have
	\begin{align}
		2\Phi^ih_{00}^{\frac{m+1}{2}}W_0^m &= h_{00}^{m+1}\{-2\sigma_0s^i_0 - \sigma_0 W_0 \bar{k}^i + \sigma_0k_0W^i + 2\sigma_0\sigma_1s_0b^i + \sigma_0\sigma_1W_r\bar{k}^rW_0b^i - \nonumber \\
		&
		\sigma_0\sigma_1k_0b^i\} +  h_{00}^{\frac{m+1}{2}}W_0^{m}{-k_0y^i + \frac{1}{2}h_{00}\bar{k}^i + \sigma_1 R_{00}b^i} + 	W_0^{2m}(-4\sigma_1 R_{00}y^i \nonumber \\
		&
		 + 2\sigma_1h_{00}W_r\bar{k}^ry^i + {\frac{\sigma_1}{2}}h_{00}^{\frac{m+3}{2}}W_r\bar{k}^rb^i),
		\end{align}
		where
		\begin{align*}
		\sigma_0 =& \frac{m}{m+1}\frac{e^{\frac{m-1}{m}k(x)}}{2^{m-1}},\\
		\sigma_1 =& \frac{2ms}{s^2 - ms^2 + mb^2}.
		\end{align*}
		Equation (3.9) can be rewritten as -
		\begin{align} \label{spray10}
		2\Phi^ih_{00}^{\frac{m+1}{2}}W_0^m = A^i_{(1)} h_{00}^{m+1} + A^i_{(2)}h_{00}^{\frac{m+1}{2}}W_0^m + A^i_{(3)}W_0^{2m},
		\end{align}
		where
		\begin{align*}
		A^i_{(1)} = & -2\sigma_0S^i_0 - \sigma_0W_0\bar{k}^i + \sigma_0k_0W^i + 2\sigma_0\sigma_1S_0W^i, \\
		A^i_{(2)} = & -k_0y^i + \frac{1}{2}h_{00}\bar{k}^i + \sigma_1R_{00}b^i,\\
		A^i_{(3)} =& -4\sigma_1R_{00}y^i + 2\sigma_1h_{00}W_r\bar{k}^ry^i + \frac{\sigma_1}{2}h_{00}^{\frac{m+3}{2}}W_r\bar{k}^rb^i.
		\end{align*}
		\section{The necessary and sufficient conditions for constant curvature of generalized Kropina spaces}
		In this section, we consider a Generalized Kropina space (M, $\alpha^{m+1}/\beta^m$) of constant curvature K, where $\alpha = \sqrt{a_{ij}y^iy^j}$ and $\beta = b_i(x)y^i$. Furthermore, we suppose that the matrix $(a_{ij})$ is always positive definite and that the dimension $n\geq2$. Hence, it follows that $\alpha^{m+1}$ is not divisible by $\beta^m$. This is an important relation and it is equivalent to that $h_{00}^{\frac{m+1}{2}}$ is not divisible by $W_0^m$ Using these, we shall obtain the necessary and sufficient conditions for a Kropina space to be of constant curvature.
		\subsection{The curvature tensor of a generalized Kropina space }
		Let ${{R_j}^i}_{kl}$ be three h-curvature tensors of Cartan connection in Finsler space. The Berwald spray curvature tensor is 
		\begin{align}
		^{(b)}{{{R_j}^i}_{kl}} = A_{(kl)} \Bigl({\frac{\partial {{G_j}^i}_k}{\partial x^l}}\ + {{G_j}^r}_k {{G_r}^i}_l\Bigr),
		\end{align}
		where the symbol $A_{(kl)}$ denotes the interchange of indices k and l and subtraction. It is well known that the equality ${{R_0}^i}_{kl} = ^{(b)}{{{R_0}^i}_{kl}}$ holds good \cite{yoshi12}.
		\par
		From 2$G^i$ =$ ^{h}{{{\gamma_0}^i}_0} + 2\Phi^i$, it follows $\quad G^i_j = ^{h}{{{\gamma_0}^i}_j} + {{\Phi_j}^i}_k$, where $\Phi^i_j$ := $\frac{\partial \Phi^i}{\partial y^j}$ and ${\Phi^i}_{jk} := \frac{\partial \Phi^i_j}{\partial y^k}$. Substituting the above equalities in (4.1), we get 
	 	 \\ $^{(b)}{{{R_j}^i}_{kl}} = ^{h}{{{R_j}^i}_{kl}} + A_{(kl)}\{{{\Phi_j}^i}_{k \lVert l} + {{\Phi_j}^r}_k {{\Phi_r}^i}_l\}$.
		
		\par
		The following results are well known \cite{Yoshi11}:
		\begin{proposition}
		The necessary and sufficient condition for a Finsler space (M, F) to be of scalar curvature K is that the equality
		\begin{align}
		{{{R_0}^i}_{0l}} = KF^2(\delta^i_l - l^il_l),
		\end{align}
		where $l^i = y^i/F$ and $l_l = \partial F/\partial y^l$, holds.
		\par
		If the equality (4.2) holds and K is constant, then the Finsler space is said to be of constant curvature K.
		\end{proposition}
	\par
	For a generalized Kropina space of constant curvature K, since F = $\varepsilon h_{00}^{\frac{m+1}{2}}/2W_0^m$, where $\varepsilon = (e^{k(x)})^{\frac{m-1}{2}}$,
	we have\\
	\\
	$l^il_l = \frac{\frac{m+1}{2}W_0 h_{0l} - mh_{00}W_l}{W_0h_{00}}y^i$.\\
	So, 
	
	$\delta^i_l - l^il_l = \delta^i_l - \frac{\frac{m+1}{2}W_0 h_{0l} - mh_{00}W_l}{W_0h_{00}}y^i$.
	\par
	Using the curvature obtained above, we have ${{R_0}^i}_{0l} = ^{h}{{{R_0}^i}_{0l}} + 2\Phi^i_{\parallel l} - \Phi^i_{l \parallel 0} + 2\Phi^r {{\Phi_r}^i}_l - \Phi^r_l\Phi^i_r$.\\
	Substituting the above equalities in (4.2), we get
	\begin{align}
	K\dfrac{\varepsilon^2h_{00}^{m+1}}{(2W_0)^{2m}}{h^i}_l = ^{h}{{{R_0}^i}_{0l}} + 2\Phi^i_{\parallel l} - \Phi^i_{l \parallel 0} + 2\Phi^r {{\Phi_r}^i}_l - \Phi^r_l\Phi^i_r.
	\end{align}
	
	\subsection{Rewriting the equation (4.3) using $h_{00}$ and $W_{0}$}
	
	\begin{enumerate}
		
	\item [(1.)] The calculation for ${\Phi^i}_{\parallel l}$.
	\par
	First, applying the h-covariant derivative $_{\parallel l}$ to (3.11), it follows:
	\begin{align}
	2h_{00}^{\frac{m+1}{2}}W_0^m\Phi^i_{\parallel l} &=& h_{00}^{m+1}A^i_{(1)\parallel l} - m h_{00}^{m+1}W_{0 \parallel l} A^i_{(1)} + (h_{00})^{\frac{m+1}{2}}W_0^m A^i_{(2)\parallel l} + \nonumber \\
	&&
	 W_0^{2m}A^i_{(3)\parallel l} + m W_0^{2m-1}W_{0 \parallel l}A^i_{(3)}. \nonumber
    \end{align}
    By appropriate substitutions, we get -
    \begin{align}
    2h_{00}^{\frac{m+1}{2}} W_{0}^m \Phi^{i}_{\parallel l} &=& h_{00}^{m+1} B^{i}_{(1) \parallel l} + mh_{00}^{m+1} B^i_{(21) l} +h_{00}^{\frac{m+1}{2}}W_{0}^m B^{i}_{(22)l} + \nonumber \\
    &&    
    W_{0}^{2m} B^i_{(3)l} + m W_0^{2m-1} B^i_{(4)l},
    \end{align}
    where
     \begin{align*}
     B^{i}_{(1) \parallel l} =& A^i_{(1)\parallel l}, \\
     B^i_{(21) l} =& -W_{0 \parallel l} A^i_{(1)}, \\
     B^{i}_{(22)l} =& A^i_{(2)\parallel l},\\
       B^i_{(3)l} =& A^i_{(3)\parallel l},\\
       B^i_{(4)l} =& W_{0 \parallel l}A^i_{(3)}.
    \end{align*}
    \item [(2.)] The calculation for $\Phi^i_l$.
    \par
    Secondly, differentiating equation (\ref{spray10}) by $y^l$, we get
%
    \begin{align}
    2\Phi^i_lh_{00}^{\frac{m+3}{2}} W_0^{m+1} &=& h_{00}^{m+2}W_0C^i_{(0)l} + h_{00}^{m+2} C^i_{(11)l} + h_{00}^{\frac{m+3}{2}} W_0^{m+1} Ci_{(12)l} + \nonumber \\
    &&
    (h_{00})^{m+1}W_0C^i_{(21)l} + h_{00}W_0^{2m+1} C^i_{(22)l} + W_0^{2m} h_{00}C^i_{(3)l} + \nonumber \\
    &&
     W_0^{2m+1}C^i_{(4)l},
    \end{align}
    where
    \begin{align*}
    C^i_{(0)l} &= A^i_{(1)l}, \\
    C^i_{(11)l} &= -mW_l A^i_{(1)}, \\
    C^i_{(12)l} &= A^i_{(2)l}, \\
    C^i_{(21)l} &= (m+1)h_{0l}A^i_{(1)}, \\
    C^i_{(22)l} &= A^i_{(3)l}, \\
    C^i_{(3)l} &= mW_lA^i_{(3)}, \\
    C^i_{(4)l} &= -(m+1)h_{0l}A^i_{(3)}.
    \end{align*}
    \item [(3).] The Calculation for $\Phi^i_{l \parallel 0}$.
    \\
Applying the h-covariant derivative $_{\parallel 0}$ to (4.5), we get
\begin{align}
2h_{00}^{\frac{m+3}{2}}W_0^{m+2} \Phi^i_{l \parallel 0} &=& h_{00}^{m+2}W_0^2 D^i_{(1)l} + h_{00}^{m+2}W_0 D^i_{(21)l} + \nonumber \\
&&
 h_{00}^{m+2} D^i_{(31)l} + h_{00}^{\frac{m+3}{2}}W_0^{m+2}D^i_{(22)l} \nonumber \\
 &&
 + h_{00}^{m+1}W_0^2D^i_{(32)l} + h_{00}^{m+1}W_0D^i_{(41)l} \nonumber \\
 &&
  + h_{00} W_0^{2m+2}D^i_{(33)l} + h_{00}W_0^{2m+1}D^i_{(42)l} + \nonumber \\
  &&
  W_0^{2m+2} D^i_{(5)l} + W_0^{2m+1}D^i_{(6)l},
\end{align}
where
\begin{align*}
D^i_{(1)l} &= C^i_{(0)l \parallel 0}, \\
D^i_{(21)l} &= C^i_{(11)l \parallel 0} - W_{0 \parallel 0}C^i_{(0)l}, \\
D^i_{(31)l} &= -2W_{0 \parallel 0}C^i_{(11)l}, \\
D^i_{(22)l} &= C^i_{(12)l \parallel 0}, D^i_{(32)l} = C^i_{(21)l \parallel 0}, \\
D^i_{(41)l} & = - W_{0 \parallel 0}C^i_{(21)}, D^i_{(33)l} = C^i_{(22)l \parallel 0}, \\
D^i_{(42)l} & = W_{0 \parallel 0}C^i_{(22)l} + C^i_{(3)l \parallel 0}, D^i_{(5)l} = C^i_{(4)l \parallel 0}, D^i_{(6)l} = W_{0 \parallel 0} C^i_{(4)l}.  
\end{align*}
\item [(4.)] The Calculation for $\Phi_l^r \Phi_r^i$
\begin{align}
& 4\Phi_l^r \Phi_r^i (h_{00})^{m+3} W_0^{2m+2} \nonumber\\
& = (h_{00})^{2m + 4} W_0^2 E^i_{(01)l} + (h_{00})^{2m+4} W_0 E^i_{(11)l} + (h_{00})^{2m+4} E^i_{(21)l} + \nonumber \\ 
&
+ (h_{00})^{\frac{3m+7}{2}}W_0^{m+2}E^i_{(12)l} +
 (h_{00})^{\frac{3m+7}{2}}W_0^{m+1}E^i_{(22)l} + (h_{00})^{2m+3} W_0 E^i_{(31)l}  \nonumber \\
 &
 + (h_{00})^{m+3}W_0^{2m+2}E^i_{(23)l} + (h_{00})^{m+3}W_0^{2m+1}E^i_{(32)l} + (h_{00})^{m+3}W_0^{2m}E^i_{(41)l} \nonumber \\
  &
+ (h_{00}^{\frac{m+5}{2}})W_0^{3m+2}E^i_{(33)l} + (h_{00})^{m+2}W_0^{2m+2}E^i_{(42)l} + 
  (h_{00})^{m+2}W_0^{2m+1}E^i_{(51)l} \nonumber \\
  &
  + (h_{00})^2W_0^{4m+2}E^i_{(43)l} +  (h_{00})^2W_0^{4m+2}E^i_{(43)l} + (h_{00})^{m+1}W_0^{2m+3}E^i_{(52)l} \nonumber \\
  &
  + (h_{00})^2W_0^{4m}E^i_{(61)l}   +   (h_{00})W_0^{4m+2}E^i_{(62)l} + (h_{00})W_0^{4m+1}E^i_{(7)l} + W_0^{4m+2}E^i_{(8)l},
\end{align}
where
\begin{align*}
E^i_{(0)l} &= C^i_{(0)r}C^r_{(0)l}, E^i_{(11)l} = C^i_{(11)r}C^r_{(0)l} + C^i_{(0)r}C^r_{(11)l}, E^i_{(21)l} = C^i_{(11)r}C^r_{(11)l}, \\
E^i_{(12)l} &= C^i_{(0)r}C^r_{(12)l} + C^i_{(12)r}C^r_{(0)l}, \\
E^i_{(22)l} &= C^i_{(12)r}C^r_{(11)l} + C^i_{(11)r}C^r_{(12)l} + C^i_{(21)r}C^r_{(0)l} + C^i_{(0)r}C^r_{(21)l},\\
E^i_{(31)l} &= C^i_{(21)r}C^r_{(11)l} + C^i_{(11)r}C^r_{(21)l}, E^i_{(23)l} = C^i_{(12)r}C^r_{(12)l} + C^i_{(22)r}C^r_{(0)l} + C^i_{(0)r}C^r_{(22)l}, \\
E^i_{(32)l} &= C^i_{(21)r}C^r_{(12)l} + C^i_{(12)r}C^r_{(21)l} + C^i_{(3)r}C^r_{(0)l} + C^i_{(22)r}C^r_{(11)l} + C^i_{(11)r}C^r_{(22)l} + C^i_{(0)r}C^r_{(3)l},\\
E^i_{(41)l} &= C^i_{(3)r}C^r_{(11)l} + C^i_{(11)r}C^r_{(3)l} + C^i_{(21)r}C^r_{(21)l}, E^i_{(33)l} =  C^i_{(22)r}C^r_{(12)l} + C^i_{(12)r}C^r_{(22)l},\\
E^i_{(42)l} &= C^i_{(4)r}C^r_{(0)l} + C^i_{(3)r}C^r_{(12)l} + C^i_{(22)r}C^r_{(21)l} + C^i_{(21)r}C^r_{(22)l} + C^i_{(12)r}C^r_{(3)l} + C^i_{(0)r}C^r_{(4)l},\\
E^i_{(51)l} &= C^i_{(3)r}C^r_{(21)l} + C^i_{(21)r}C^r_{(3)l} + C^i_{(4)r}C^r_{(11)l} + C^i_{(11)r}C^r_{(4)l},
E^i_{(43)l} = C^i_{(22)r}C^r_{(22)l},\\
E^i_{(52)l} &= C^i_{(4)r}C^r_{(12)l} + C^i_{(12)r}C^r_{(4)l} + C^i_{(3)r}C^r_{(22)l} + C^i_{(22)r}C^r_{(3)l},\\
E^i_{(61)l} &= C^i_{(3)r}C^r_{(3)l} + C^i_{(4)r}C^r_{(21)l} + C^i_{(21)r}C^r_{(4)l}, E^i_{(62)l} = C^i_{(4)r}C^r_{(22)l} + C^i_{(22)r}C^r_{(4)l},\\
E^i_{(7)l} &= C^i_{(4)r}C^r_{(3)l} + C^i_{(3)r}C^r_{(4)l}, E^i_{(8)l} = C^i_{(4)r}C^r_{(4)l}.
\end{align*}
\item [(5.)] The Calculation for $\Phi^r {{\Phi_r}^i}_l$.
\\
Differentiating (4.5) by $y^r$, we get
\begin{align}
&4h_{00}^{m+3}W_0^{2m+2}\Phi^r {{\Phi_r}^i}_l = h_{00}^{2m+4}W_0J^i_{(11)l} + h_{00}^{2m+4}J^i_{(21)l} + \nonumber \\
&
 h_{00}^{2m+3}W_0^{m+2}J^i_{(12)l} + h_{00}^{2m+3}W_0^{m+1}J^i_{(22)l} + 
 h_{00}^{2m+3}W_0J^i_{(31)l} + \nonumber \\
 &
 h_{00}^{2m+3}W_0^{2m+2}J^i_{(23)l}  
  + h_{00}^{\frac{3m+5}{2}}W_0^{2m+1}J^i_{(32)l} +\nonumber \\
  &
   h_{00}^{m+3}W_0^{m+1}J^i_{(41)l} + 
   h_{00}^{\frac{3m+5}{2}}W_0^{3m+2}J^i_{(33)l} + h_{00}^{m+2}W_0^{2m+2}J^i_{(42)l} \nonumber \\
   &
   + h_{00}^{m+2}W_0^{m+2}J^i_{(51)l} + h_{00}^{m+1}W_0^{4m+2}J^i_{(43)l} 
   + h_{00}^{\frac{m+3}{2}}W_0^{3m+2}J^i_{(52)l} +
   \nonumber \\
   & h_{00}^{\frac{m+3}{2}}W_0^{2m+2}J^i_{(61)l} + 
    h_{00}W_0^{3m+2}J^i_{(71)l} + W_0^{4m+2}J^i_{(8)l}.
\end{align}
\item [(6.)] The main relation\\
Using the equation (4.3)
\begin{align*}
	K\dfrac{\varepsilon^2h_{00}^{m+1}}{(2W_0)^{2m}}{h^i}_l = ^{h}{{{R_0}^i}_{0l}} + 2\Phi^i_{\parallel l} - \Phi^i_{l \parallel 0} + 2\Phi^r {{\Phi_r}^i}_l - \Phi^r_l\Phi^i_r.
	\end{align*}
	Multiplying by $h_{00}^{m+2}W_0^{2m+2}$, we have the equality
	\begin{align*}
	& 2^{2m}Kh_{00}^{2m+4}W_0^{m+1}h^i_l = 2^{2m+2}h_{00}^{2m+2}W_0^{2m+2}  + \nonumber \\
	&
	^{h}{{{R_0}^i}_{0l}} 2^{2m+1}h_{00}^{2m+1}W_0^{m+1}.2^{2m}W_0^{m+1}
	 \Phi^i_{\parallel l} - 
	 2^{2m}h_{00}^{m+1}W_0.2^{2m}h_{00}^{2m}W_0^{2m+1}\Phi^i_{l \parallel 0} + \nonumber \\
	 &
	  2^{2m+3}h_{00}^{2m+2}W_0^{2m+2}\Phi^r{{{\Phi_r}^i}_l},
	\end{align*}
	where ${h^i}_l = {\delta^i}_l - l^il_l$. Using $\Phi^i_{\parallel l}, \Phi^i_l, \Phi^i_{l \parallel 0}, \Phi^r{{\Phi_r}^i}_l, {\Phi^r}_l{\Phi^i}_r$ in the above equality, by straight forward computation, we finally obtain
	\begin{align}
	h_{00}^{2m+2}{{\gamma_1}^i}_{(2m+3)l} + h_{00}^{m+1}{{\gamma_2}^i}_{(2m+7)l} + W_0^{2m+2}{{\gamma_3}^i}_{(2m+7)l} = 0,
	\end{align} 
	where ${{\gamma_1}^i}_{(2m+3)l}, {{\gamma_2}^i}_{(2m+7)l}$ and ${{\gamma_3}^i}_{(2m+7)l}$ are homogeneous polynomials of degree 2m+3, 2m+7 and 2m+7 in $y^i$ respectively. They are called the curvature part, the vanishing part and the killing part, respectively.
	
\end{enumerate}

	\begin{proposition}
	The necessary and sufficient condition for a Kropina space (M, F) with $F = \frac{\alpha^{m+1}}{\beta^m} = \frac{(e^{k(x)})^{\frac{m-1}{2}}h_{00}^{\frac{m+1}{2}}}{2^m {W_0}^m}$
	\end{proposition}
	to be of constant curvature K is that (4.9) holds good.
	
	\end{document}